\documentclass[a4paper]{article}

\usepackage{epsfig}

\def\refeq#1{(\ref{eq:#1})}

\newcommand{\id}{{\rm Id}}

\newcommand{\grad}{{\rm grad\;}}

\newcommand{\R}{{\bf R}}

\newcommand{\N}{{\bf N}}

\newtheorem{thm}{Theorem}[section]

\newtheorem{lem}{Lemma}[section]
\newtheorem{cor}{Corollary}[section]
\newtheorem{definition}{Definition}[section]
\newtheorem{rem}{Remark}[section]

\def\thebibliography#1{\section*{Reference}\list
 {[\arabic{enumi}]}{\settowidth\labelwidth{[#1]}\leftmargin\labelwidth
 \advance\leftmargin\labelsep
 \usecounter{enumi}}
 \def\newblock{\hskip .11em plus .33em minus .07em}
 \sloppy\clubpenalty4000\widowpenalty4000
 \sfcode`\.=1000\relax}

\def\resume{\if@twocolumn
\section*{Abstract}
\else \small 
\begin{center}
{\bf R\'esum\'e\vspace{-.5em}\vspace{0pt}} 
\end{center}
\quotation 
\fi}
\def\endresume{\if@twocolumn\else\endquotation\fi}

\newenvironment{thmintro}{\par\medskip\noindent{\bf 
Theorem}\begingroup\em}{\endgroup\hfill\par\medskip} 
\newenvironment{corintro}{\par\medskip\noindent{\bf 
Corollary}\begingroup\em}{\endgroup\hfill\par\medskip} 
\newenvironment{preuve}{\par\noindent{\bf Proof~:}}{\hfill\stopthm} 
\newcommand{\stopthm}{\hfill$\diamond$\par\smallskip} 
\author{St\'ephane Grognet}
\title{Marked length spectrum of magnetized surfaces}
\begin{document}

\maketitle

\vskip 30 pt
\noindent Universit\'e de Nantes, D\'epartement de Math\'ematiques, Laboratoire Jean Leray U. M. R. 6629, 2, rue
de la Houssini\`ere, BP 92208, F-44322 Nantes cedex 03.
 
\noindent Stephane.Grognet@univ-nantes.fr

\begin{abstract}
The main result presented here is that the flow associated with 
a riemannian metric and a non zero magnetic field 
on a compact oriented surface without boundary,
under assumptions of hyperbolic type,
cannot have 
the same length spectrum of topologically corresponding periodic orbits as
the geodesic flow associated with another riemannian metric  
having a negative curvature and the same total volume.
The main tool is a regularization inspired by U. Hamenst\"adt's methods.
\end{abstract}
%\tableofcontents 

%%%%%%%%%%%%%%%%%%%%%%%%%%%%%%%%%%%%%%%%%

\section{Introduction}%\label{sec:}

%%%%%%%%%%%%%%%%%%%%%%%%%%%%%%%%%%%%%%%%%

The problems of entropic and spectral rigidity of riemannian manifolds
have been widely studied, beginning with the surfaces~\cite{Ka3}. 
The works treat riemannian metrics on compact surfaces~\cite{Cr}, \cite{Cr-Fa-Fe}, \cite{Ot1}, 
on higher dimension manifolds~\cite{Be-Co-Ga1}, \cite{Be-Co-Ga2}, \cite{Cr-Kl}, \cite{Fa}, 
or on surfaces with singularities~\cite{He-Pa}.
The related problem of boundary rigidity of a riemannian metric features many results~\cite{L-S-U},~\cite{S1},~\cite{S2}.
The rigidity of an absolutely continuous flow conjugacy persists in some way
with the presence of a magnetic field on a compact surface~\cite{Gro1},
and so do entropic rigidity in this case~\cite{Gro2}.
The topological entropy of the magnetic flow in higher dimension has also been studied~\cite{Pa3},~\cite{B-P}.

Unlike the geodesic flow, a conjugacy being only continuous (in fact H\"older-continuous) between two magnetic flows 
on a surface had not been treated.

When the surface is compact and the Jacobi endomorphism~\cite{Fo1} of the magnetic flow is negative,
this flow has got the Anosov property~\cite{Go}~; two such flows have got the same marked length spectrum
of periodic orbits if and only if they are~${\cal C}^0$-conjugated~\cite{Gro1}.

The main result presented here is that the flow associated with a riemannian metric and a non zero magnetic field on a closed surface,
if it has got a negative Jacobi endomorphism,
cannot have the same marked length spectrum as the geodesic flow associated with another riemannian metric  
having a negative curvature and the same total volume.
The assumption on the equality of the total volumes is essential~\cite{Gro1}.

\begin{thmintro} {\bf \ref{thm:spm}}
Let~$M$ be a closed (compact without boundary), connected, oriented surface.
Let~$g_1$ and~$g_2$ be two~${\cal C}^\infty $-riemannian metrics over~$M$
whose curvatures are negatively pinched~:~$-k_0^2\leq K_i\leq -k_1^2<0$ for~$i=1,2$,.
Let~$\kappa _1$ be a~${\cal C}^\infty $-magnetic field over~$M$.
The magnetic flow~$\psi ^1_t=\psi ^{g_1,\kappa _1}_t$ is supposed to have a negative Jacobi endomorphism.
If the magnetic flow~$\psi ^1_t$ and the geodesic flow~$\varphi ^2_t=\varphi ^{g_2}_t$ have the same marked length spectrum,
and if the surface~$M$ has the same total volume for the two metrics,
then the two metrics are isotopic, which means that one is the image of the other by a diffeomorphism~$f$ of~$M$
homotopic to the identity, and the magnetic field~$\kappa _1$ is zero.
\end{thmintro}

The proof consists in coming back to the known case where there exists an absolutely continuous conjugacy between the two flows~\cite{Gro1}.
The proof of the regularity of the conjugacy is inspired by U. Hamenst\"adt's methods~\cite{Ha1,Ha2}. 
We construct linearizations of the universal covering of the surface, compatible with the stable spaces of the flow.
This is useful to proof that the Lyapounoff exponents of the periodic orbits are preserved (theorem~\ref{thm:reguconj}), 
which ensures that the conjugacy is smooth~\cite{L-M}. 
The regularity of the conjugacy used to proof the theorem~\ref{thm:spm} is valid in general for two magnetic flows (we denote~$T^1_iM$ the unit tangent bundle of~$g_i$)~:

\begin{corintro} {\bf \ref{cor:reguconj}}
Let~$M$ be a closed (compact without boundary), connected, oriented surface.
Let~$g_1$ and~$g_2$ be two~${\cal C}^\infty $-riemannian metrics over~$M$
whose curvatures are negatively pinched~:~$-k_0^2\leq K_i\leq -k_1^2<0$ for~$i=1,2$.
Let~$\kappa _1$,~$\kappa _2$ be two~${\cal C}^\infty $-magnetic fields over~$M$.
The two magnetic flows~$\psi ^1_t=\psi ^{g_1,\kappa _1}_t$ and~$\psi ^2_t=\psi ^{g_2,\kappa_2}_t$ 
are supposed to have negative Jacobi endomorphisms.
If the two magnetic flows have the same marked length spectrum,
then they are conjugated by a~${\cal C}^\infty $-diffeomorphism~$h$ from~$T^1_1M$ onto~$T^1_2M$.
\end{corintro}

A uniformization of a surface equipped with a metric with negative curvature has already been constructed~\cite{Fe-Or}~;
it applies to an Anosov flow on a~$3$-manifold, but with the condition that the stable spaces be of~${\cal C}^1$-class,
which is unlikely for the magnetic flow~\cite{Pa2}.
The uniformly quasiconformal diffeomorphisms present another example of uniform structures on stable spaces~\cite{K-S}.

It seems legitimate to ask if the construction presented here is practicable for other flows
whose stable spaces are not necessarily of~${\cal C}^1$-class,
and particularly to which extent a~${\cal C}^0$-conjugacy between two such flows could be differentiable.

Given a manifold~$M$, diffeomorphic to~$\R^2$, with a magnetic flow having a negative Jacobi endomorphism
and having the gradient of centre-stable and centre-unstable spaces uniformly bounded,
and given a point~$p\in M$, and a unitary vector~$v\in T^1_pM$, 
the {\em linearization}~$E_v$ (defined in section~\ref{sec:linea}) sends~$M$ onto~$T_pM$~;
the geodesic directed by~$v$ onto the straight line~$\R v$~;
and the horocycles associated with the centre-stable manifold of~$v$ onto the straight lines orthogonal to~$v$.
This linearization, used as is, presents a little rigidity.

\begin{thmintro} {\bf \ref{thm:rigiR}}
Let~$M$ be an oriented surface diffeomorphic to~$\R^2$, equipped with two~${\cal C}^\infty $-riemannian metrics~$g_1$ and~$g_2$
whose curvatures are negatively pinched~:~$-k_0^2\leq K_i\leq -k_1^2<0$ for~$i=1,2$. 
Let~$\kappa _1$,~$\kappa _2$ be two~${\cal C}^\infty $-magnetic fields over~$M$.
The two magnetic flows~$\psi ^1_t=\psi ^{g_1,\kappa _1}_t$ and~$\psi ^2_t=\psi ^{g_2,k_2}_t$ 
are supposed to have pinched negative Jacobi endomorphisms
(the~${\cal C}^1$-norms of~$\kappa _1$ and~$\kappa _2$ are thus bounded),
and the gradient of the centre-stable~$u_{-,i}$ spaces and the gradient of the centre-unstable~$u_{+,i}$ spaces for~$i=1,2$
are supposed to be uniformly bounded.
If there exist a diffeomorphism~$f:M\to M$ and a point~$p\in M$ satisfying
\[
\forall v \in T^1_{1,p}M \qquad
\exists v' \in T^1_{2,f(p)}M \qquad
E^2_{v'} \circ f = E^1_v ,
\]
then the two metrics are images one of each other by~$f$, and so are the two magnetic fields~:~$\kappa_2=\kappa_1 \circ f$.
\end{thmintro}

When metrics and magnetic fields are invariants under a cocompact group, the rigidity of the linearization is stronger in some way.

\begin{thmintro} {\bf \ref{thm:rigicomp}}
Let~$M$ be a closed (compact without boundary), connected, oriented surface.
Let~$g_1$ and~$g_2$ be two~${\cal C}^\infty $-riemannian metrics over~$M$
whose curvatures are negatively pinched~:~$-k_0^2\leq K_i\leq -k_1^2<0$ for~$i=1,2$.
Let~$\kappa _1$,~$\kappa _2$ be two~${\cal C}^\infty $-magnetic fields over~$M$.
The two magnetic flows~$\psi ^1_t=\psi ^{g_1,\kappa _1}_t$ and~$\psi ^2_t=\psi ^{g_2,\kappa _2}_t$ 
are supposed to have negative Jacobi endomorphisms.
If there exist two vectors~$v_1\in T^1_1\widetilde M$,~$v_2\in T^1_2\widetilde M$ 
and a~${\cal C}^1$-diffeomorphism~$f:M\to M$ homotopic to the identity,
of which a lift~$\widetilde f$ over~$\widetilde M$ satisfies~$E^2_{v_2}\circ \widetilde f = E^1_{v_1}$,
then the two metrics are isotopic, transported by~$f$,
and so are the two magnetic fields~:~$\kappa_2=\kappa_1 \circ f$.
\end{thmintro}

Using the tools of the construction of the linearization,
we also get a result of constancy of (future) Lyapounoff exponents along the centre-stable manifolds (theorem~\ref{thm:Lyap}).

%%%%%%%%%%%%%%%%%%%%%%%%%%%%%%%%%%%%%%%%%%%%%%%%%%%

\section{Notations and background}\label{sec:notat}

%%%%%%%%%%%%%%%%%%%%%%%%%%%%%%%%%%%%%%%%%%%%%%%%%%%

In the following,~$M$ denotes a complete, connected, oriented surface,
equipped with a~${\cal C}^\infty $-riemannian metric 
whose curvature is negatively pinched~:~$-k_0^2\leq K\leq -k_1^2<0$.
The Cartan-Hadamard theorem~(\cite{Ga-Hu-La}, p.138) implies that 
the universal cover~$\widetilde{M}$ is diffeomorphic to~$\R ^2$, with cover mapping~$\Pi : \widetilde M \to M$.
Within sections~\ref{sec:courbure},~\ref{sec:fluc},~\ref{sec:transport},~\ref{sec:linea} and~~\ref{sec:regu},
the surface~$M$ is simply connected, thus equal to~$\widetilde M$.
The projection of~$TM$ and~$T^1M$ on~$M$ is written down~$\pi$. 
The cover mapping is~$\Pi : \widetilde M \to M$.

The  surface~$M$ is said {\em closed} if it is compact (without boundary).

Let~$N$ be the rotation of angle~$+\pi /2$ in the tangent space~$TM$.

For a curve~$c:\R \to M$, the equation of the magnetic flow~$\psi_t=\psi^{g,\kappa}_t=\psi^{\kappa}_t$ associated with a magnetic field~$\kappa :M\to \R $ is~\cite{Gro1}~:
\[
\frac{Dc}{dt} = \kappa (c(t)) \ N\left( \frac{dc}{dt}\right) .
\]
The flow is a one-parameter group of diffeomorphisms acting on~$T^1M$. The magnetic field~$\kappa $ is supposed to be smooth.

The {\em Jacobi endomorphism} associated with this second order differential equation~\cite{Fo1} is the application~\cite{Gro1}~:
\[
\matrix{
q & : & T^1M & \to & \R \cr
{} & {} & v & \mapsto & K(\pi (v)) +\kappa (\pi (v))^2 -\left\langle N(v) , (\grad \kappa ) (\pi (v)) \right\rangle .
}
\]
When the surface~$M$ is compact (closed), 
saying that the Jacobi endomorphism is negative is equivalent to
saying that it is pinched between two strictly negative constants.
In the following, the real function~$\kappa :M\to \R $ is a magnetic field 
such that the associated Jacobi endomorphism~$q$ satisfies the pinching condition,
which means that there exist two positive constants~$q_0$ et~$q_1$ verifying~:
\begin{equation}\label{eq:pince}
-q_0^2\leq q\leq -q_1^2 <0 .
\end{equation}

\begin{definition}~{\bf \cite{Gro1}}
With the assumption~\refeq{pince}, to a vector~$v\in T^1M$ are associated
the stable~$j_-(v,t)$ and unstable~$j_+(v,t)$ Jacobi fields along the orbit of~$v$,
with components~$(x_-(v,t),y_-(v,t))$ and~$(x_+(v,t), y_+(v,t))$ in the base~$(\psi_tv,N(\psi_tv))$
satisfying
\[
\matrix{
x_-(v,+\infty )=0, & y_-(v,+\infty )=0, & y_-(v,0)=1\, ; \cr
x_+(v,-\infty )=0, & y_+(v,-\infty )=0, & y_+(v,0)=1. \cr
}
\]
The stable and unstable spaces are determined by the mappings~:
\[
v\mapsto \left( w_-(v),u_-(v) \right) = \left( x_-(v,0),{y_-}'(v,0) \right) , 
\]
\[
v\mapsto \left( w_+(v),u_+(v) \right) = \left( x_+(v,0),{y_+}'(v,0) \right) . 
\]
\end{definition}

The tangential component of the stable space at~$v\in T^1M$ satisfies the relation~\cite{Gro1}~:
\begin{equation}\label{eq:jacosta1}
w_-(v)=x_-(v,0)=\int_{t=+\infty}^0 \kappa (\pi \psi ^\kappa _tv) y_-(v,t) \ dt .
\end{equation}

Writing
\begin{equation}\label{eq:C1}
C_1= 1+\frac{\Vert \kappa \Vert _\infty +\Vert \kappa \Vert _\infty ^2}{q_1} +q_0 ,
\end{equation}
yields as in~\cite{Gro2},~section~~$3.2$~:
\begin{equation}\label{eq:jacosta2}
\Vert j_-(v,t)\Vert + \Vert j_-'(v,t)\Vert \leq C_1 y_-(v,t) .
\end{equation}

Let~$W^S(v)$,~$W^{CS}(v)$,~$W^U(v)$ and~$W^{CU}(v)$ be respectively 
the stable, centre-stable, unstable, centre-unstable manifolds associated with the unitary vector~$v$.
The stable horocycle of the magnetic flow associated with~$v$ is~$H^\kappa _v(0)=H_v(0)=\pi W^S(v)$.
The stable horocycle associated with~$\psi_tv$ is~$H^\kappa _v(t)=H_v(t)$.
The Busemann function associated with~$v$ is the mapping~$B_v:M\to \R $ such that~$B_v(H_v(t))=t$.
Under the assumption~\refeq{pince}, the centre-stable and centre-instable spaces identified to
the Ricatti applications~$u_-(v)=y_-'(v,0)$ and~$u_+(v)=y_+'(v,0)$ are of~${\cal C}^1$-class over~$T^1M$ and~$T^1\widetilde M$~; 
the horizontal (orthogonal) component of the stable~($(v,t)\mapsto y_-(v,t)$) and instable~($(v,t)\mapsto y_+(v,t)$) jacobi fields
are of~${\cal C}^1$-class over~$T^1M\times \R $ and~$T^1\widetilde M\times \R $~; 
the circle at infinity~$\partial \widetilde M$ admits also a differential structure of~${\cal C}^1$-class~(\cite{Gro1}, section~$7$).

Let~$v^\kappa _{+\infty}=v_{+\infty}$ be the point at infinity corresponding to the future orbit of~$v \in T^1\widetilde M$ 
(and identified with its centre-stable manifold).
Given two distinct points~$x\in M$ and~$y\in M\cup \partial M$, we denote~$v^\kappa (x,y)=v(x,y)$ as the unitary vector tangent to~$M$ at~$x$, 
directing the unique curve solution joining~$x$ to~$y$ (in this order).

%%%%%%%%%%%%%%%%%%%%%%%%%%%%%%%%%%%%%%%%%%%%%%%%%%%

\section{Liouville measure and symplectic structure on the space of orbits}\label{sec:symp}

%%%%%%%%%%%%%%%%%%%%%%%%%%%%%%%%%%%%%%%%%%%%%%%%%%%

Let~$M$ be a complete connected oriented surface, equipped with a riemannian metric~$g$
of~${\cal C}^\infty $-class with pinched negative curvature~$-k_0^2\leq K\leq -k_1^2$, 
and with a uniformly bounded magnetic field~$\kappa $, of~${\cal C}^\infty $-class, 
with its Jacobi endomorphism satisfying the pinching condition~\refeq{pince}.
Using the method of the second order differential equations of Foulon~\cite{Gro2}, 
let~$X(\kappa )$ be the generating field of the magnetic flow,~$H_0$ be the horizontal field,~$Y_0$ be the vertical field~;
together, they constitute a basis field of the bundle~$TT^1M$, tangent to the unitary tangent bundle.
Let~$(X(\kappa )^\star ,H_0^\star ,Y_0^\star )$ be the dual basis field.
Let~$v$ be in $T^1M$ and~$j_i=x_i\psi ^\kappa _tv+y_iN(\psi ^\kappa _tv)$ for~$i=1,2$ be two Jacobi fields, with~$x'_i=\kappa y_i$~;
they are associated to the tangent vectors~$\xi _i=x_iX(\kappa )+y_iH_0+y'_iY_0$ at every point of the orbit of~$v$.
The wronskian form
\[
\Omega =H_0^\star \wedge Y_0^\star 
\]
applied to the pair of Jacobi fields yields
\[
\Omega (j_1,j_2)=y_1y'_2-y_2y'_1 .
\]
This is an invariant~$2$-form under the magnetic flow on~$T^1M$.
It gives a symplectic form on the space of orbits of the magnetic flow.
Its absolute value is equal to the Liouville measure on the space of orbits~\cite{Gro1,Gro2}.
We have~:
\begin{equation}\label{eq:sympl}
\Omega \left( j_-(v,\cdot ),j_+(v,\cdot ) \right) = (u_+-u_-)(v) = y_-(v,t)(u_+-u_-)(\psi ^\kappa _tv)y_+(v,t)
\end{equation}
for every~$t\in \R$. 
This quantity lays between~$2q_1$ and~$2q_0$.

%%%%%%%%%%%%%%%%%%%%%%%%%%%%%%%%%%%%%%%%%%%%%%%%%%%

\section{Curvature of the horocycles}\label{sec:courbure}

%%%%%%%%%%%%%%%%%%%%%%%%%%%%%%%%%%%%%%%%%%%%%%%%%%%

\begin{thm}\label{thm:courbure}
On a complete, connected, simply connected, oriented surface~$M$, equipped with a~${\cal C}^\infty $-riemannian metric~$g$
whose curvature is negatively pinched~: $-k_0^2\leq K\leq -k_1^2<0$, and with a~${\cal C}^\infty $-magnetic field~$\kappa $
whose~${\cal C}^1$-norm is bounded, with Jacobi endomorphism satisfying the pinching condition~\refeq{pince},
and such that the gradient of the centre-stable~$u_-$ (respectively centre-unstable~$u_+$) spaces is uniformly bounded,
the geodesic curvature of the stable (respectively unstable) horocycles of the magnetic flow
is uniformly bounded.
\end{thm}
\begin{preuve}
\begin{figure}[htb]
%\vskip 12 cm
\begin{center}
\epsfig{file=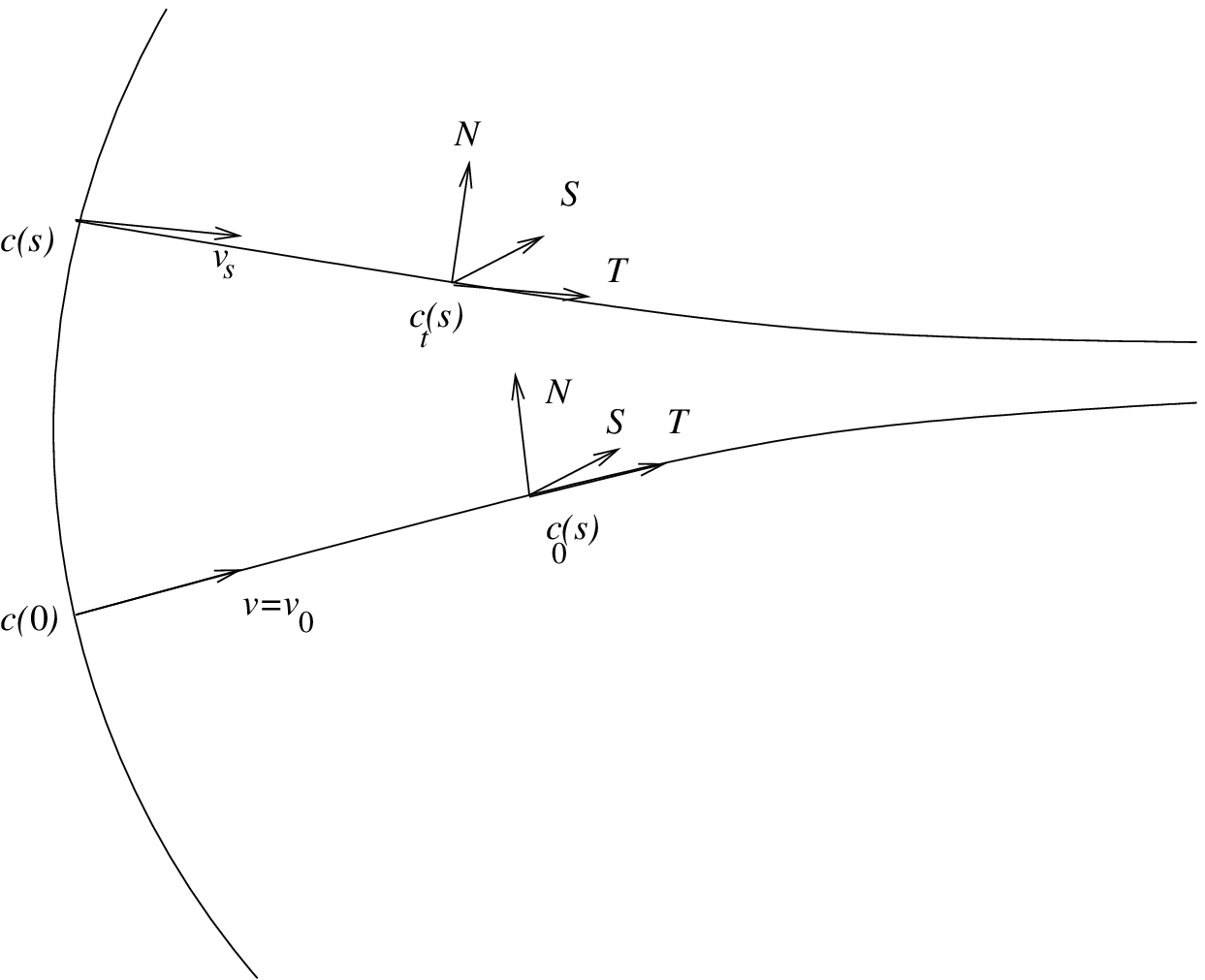}
\end{center}
\label{fig:Un}
\caption{}
\end{figure}
The unstable case is similar to the stable case~; we only consider this last one. Let~:
\begin{list}{}{}
\item{i)} $v$ be in $T^1M$,
\item{ii)} $s\to v_s$ be a smooth curve from~$]-\varepsilon ,\varepsilon [$ to~$W^{S}(v)$
such that~$\left\langle \frac{\partial}{\partial s}\pi v_s,N(v_s) \right\rangle=1$ for every~$s\in ]-\varepsilon ,\varepsilon [$,
\item{iii)} $c(s)=\pi v_s$,
\item{iv)} $c_t(s)=\pi \psi_t^\kappa v_s$,
\item{v)} $T(s,t)=\frac{D}{dt}c_t(s)$,~$N(s,t)=N(T(s,t))$,
\item{vi)}~$S(s,t)=\frac{D}{ds}c_t(s)=x_-(v,t)T+y_-(v,t)N$,
\item{vii)} ${\cal Y}(s,t)=\frac{\partial}{\partial s}y_-(v_s,t)$,
\item{viii)} $u_-(\psi_t^\kappa v_s)=y_-(v_s,t)^{-1}\frac{\partial}{\partial t} y_-(v_s,t)
=\frac{\partial}{\partial t} \ln y_-(v_s,t)$,
\item{ix)} ${\cal W}(s,t)=y_-(v_s,t)^{-1}{\cal Y}(s,t)=\frac{\partial}{\partial s} \ln y_-(v_s,t)$,
\item{x)} $L_0^s(c_t)$ be the length of the curve~$c_t$ on the interval~$[0,s]$.
\end{list}
We have~${\cal W}(s,0)=0$ for all~$s$ and
\[
{\cal W}'(s,t)=\frac{\partial}{\partial t}{\cal W}(s,t)
=\frac{\partial}{\partial s} u_-(\psi_t^\kappa v_s) .
\]
The equations~\refeq{C1},~\refeq{jacosta2}, imply
\begin{equation}\label{eq:torsion1}
|{\cal W}'(s,t)| \leq \Vert \nabla u_- \Vert_\infty \ 
C_1 \ y_-(v_s,t)  .
\end{equation}
Thus we have for all~$t\geq 0$ and all~$s$~:
\[
|{\cal W}(s,t)| \leq C_1\Vert \nabla u_- \Vert_\infty 
\int_{\tau =0}^t e^{-q_1\tau } \ d\tau  
\leq \frac{C_1\Vert \nabla u_- \Vert_\infty}{q_1}  ,
\]
\begin{equation}\label{eq:fluc1}
|{\cal Y}(s,t)|\leq 
\frac{C_1\Vert \nabla u_- \Vert_\infty}{q_1} y_-(v_s,t)  ,
\end{equation} 
With the assumption~$\left\langle S,N(T)\right\rangle (s,0)=1$ 
for all~$s\in ]-\varepsilon ,\varepsilon [$, there comes~: 
\begin{equation}\label{eq:S1}
\frac{dc}{ds}(s)=S(s,0)=w_-(v_s)v_s+N(v_s) ,
\end{equation}
thus
\[ %%begin{equation}\label{eq:S2}
1\leq \left\Vert \frac{dc}{ds}(s) \right\Vert \leq C_1 .
\] %%end{equation}
We have also~:
\begin{equation}\label{eq:S3}
\frac{Dv_s}{ds}=\left. \frac{D}{dt} \right\rfloor_{t=0} 
\hskip -15 pt 
S=(y_-'+\kappa x_-)(s,0)N(v_s)
=\left( u_-(v_s)+\kappa (c(s))w_-(v_s) \right) N(v_s) .
\end{equation}
From the relations~\refeq{S1} and~\refeq{S3} comes~:
\begin{equation}\label{eq:accel}
\frac{D}{ds} \frac{dc}{ds} (s)= \left[ \frac{dw_-(v_s)}{ds} -(u_-+\kappa w_-) \right] v_s
+ (u_-+\kappa w_-)w_- N(v_s) .
\end{equation}
In order to estimate the curvature of~$c$, the norm of~$dw_-(v_s)/ds$ should be controlled.
The equation~\refeq{jacosta1} gives~:
\[
\frac{dw_-(v_s)}{ds}=\int_{t=+\infty}^0 \left( \left\langle \nabla \kappa ,S\right\rangle  y_-(v_s,t)
+\kappa {\cal Y}(s,t) \right) \ dt .
\]
The upper estimate~\refeq{jacosta2} implies~:
\[
\left\Vert S (s,t) \right\Vert \leq 
C_1 y_-(v_s,t) \leq C_1 e^{-q_1t} .
\]
With the equation~\refeq{fluc1}, this yields for all~$s\in [-\varepsilon ,\varepsilon ]$~:
\[
\left| \frac{dw_-(v_s)}{ds} \right| \leq 
\frac{\Vert \nabla \kappa \Vert_\infty }{2q_1} C_1
+\Vert \kappa \Vert_\infty \frac{C_1 \Vert \nabla u_- \Vert_\infty }{q_1^2} .
\]
Let this last constant be written down~$C_2$. The geodesic curvature of the  horocycle~$c(s)$ is
\[
\kappa _-(s) =
\frac{\det \left( \displaystyle \frac{dc}{ds},\frac{D}{ds}\frac{dc}{ds} \right) }
{ \displaystyle \left\Vert \frac{dc}{ds}^3\right\Vert } .
\]
From the equalities~\refeq{S1} and~\refeq{accel} we deduce
\[
\kappa _- ( 1+w_-^2)^{\frac32} = (1+w_-^2)(u_-+\kappa w_-)-\frac{dw_-(v_s)}{ds} ,
\]
thus
\[
|\kappa _-| \leq 
\left( q_0 + \Vert \kappa \Vert _\infty \Vert w_- \Vert _\infty \right) 
+ C_2  .
\]
In conclusion, the geodesic curvature of the stable horocycles is uniformly bounded, 
and so is the geodesic curvature of the unstable horocycles.
\end{preuve}

The proof of the following result is left to the reader as an exercise.
\begin{cor}\label{cor:courbure}
With the assumptions of the theorem~\ref{thm:courbure}, there exists a mapping~$f:\R^+\to \R ^+$,
continuous at~$0$ and which annulates at~$0$, such that for all horocycle~$H$ of the magnetic flow,
every diffeomorphism~$c:\R \to H$ and all~$a,b\in \R $, we have
\[
L_a^b(c) \leq f\left( d \left( c(a),c(b) \right) \right) .
\]
\end{cor}

%%%%%%%%%%%%%%%%%%%%%%%%%%%%%%%%%%%%%%%%%%%%%%%%%%%

\section{Fluctuation of the stable Jacobi fields}\label{sec:fluc}

%%%%%%%%%%%%%%%%%%%%%%%%%%%%%%%%%%%%%%%%%%%%%%%%%%%

\begin{thm}\label{thm:fluc}
Let~$M$ be a complete, connected, oriented, simply connected surface, equipped with a~${\cal C}^\infty $-riemannian metric~$g$
whose curvature is negatively pinched~:~$-k_0^2\leq K\leq -k_1^2<0$, and with a~${\cal C}^\infty $-magnetic field~$\kappa $
whose~${\cal C}^1$-norm is bounded, with Jacobi endomorphism satisfying the pinching condition~\refeq{pince},
and such that the gradient of the centre-stable space~$u_-$ is uniformly bounded.
For~$\theta \in \partial M$ and~$(p',p)\in M^2$, the family of mappings~$y_-(v(p',\theta ),t)/y_-(v(p,\theta ),t)$
converges when~$t$ grows to~$+\infty $, uniformly on the compacts of~$\partial M \times M^2$, to a mapping~${\cal X}(\theta ,p,p')$,
continuous on~$\partial M\times M^2$, and which admits a partial derivative with respect to~$p,p'$ in the direction~$M^2$, 
continuous on~$\partial M\times M^2$.
\end{thm}
\begin{preuve}
Let~:
\begin{list}{}{}
\item{i)} $\theta $ belong to $\partial M$,
\item{ii)} $p$ belong to $M$,
\item{iii)} $c(s)$ be a smooth parametrization of the horocycle~$\pi W^S(v(p,\theta ))$ such that~$c(0)=p$ 
and~$\left\langle \frac{dc}{ds} , N(v(c(s,\theta )) \right\rangle =1$,
\item{iv)} $p_\theta (r,s)=\pi \psi ^\kappa _rv(c(s),\theta )$ 
the smooth parametrization of~$M$ by~$\R ^2$ which follows from it for all~$\theta \in \partial M$,
\item{v)} $T_\theta (s,t)=\frac{\displaystyle \partial }{\displaystyle \partial t}p_\theta (t,s)$,~$N_\theta (s,t)=N(T_\theta (s,t))$,~$S_\theta (s,t)
=\frac{\displaystyle \partial }{\displaystyle \partial s}p_\theta (t,s)$,
\item{vi)} ${\cal Z}(\theta ,r,s,t)=\frac{\displaystyle y_-(v(p_\theta (r,s),\theta ),t)}{\displaystyle y_-(v(p,\theta ),t)}$,
\item{vii)} ${\cal Y}(\theta ,s,t)=\frac{\displaystyle \partial}{\displaystyle \partial s}y_-(v(c(s),\theta ),t)$,
\item{viii)}  ${\cal W}(\theta ,s,t)=y_-(v(c(s),\theta ),t)^{-1}{\cal Y}(\theta ,s,t)
=\frac{\displaystyle \partial}{\displaystyle \partial s} \ln y_-(v(c(s),\theta ),t)$,
\item{ix)} $L_a^b(p_\theta (r,.))$  be the length of the curve~$s\mapsto p_\theta (r,s)$ on the interval~$[a,b]$.
\end{list}
For all real numbers~$r,s,t$ we have~:
\begin{equation}\label{eq:decalage}
{\cal Z}(\theta ,r,s,t)=\frac{y_-(v(p_\theta (r,s),\theta ),-r)}{y_-(v(p_\theta (r,0),\theta ),-r)} {\cal Z}(\theta ,0,s,t+r) .
\end{equation}
The mapping~$(\theta ,r,s,t)\mapsto y_-(v(p_\theta (r,s),\theta ),t)$ is continuous on~$\partial M\times \R ^3$ and admits partial derivatives
with respect to~$r,s,t$, continuous on~$\partial M\times \R ^3$.
Studying the uniform convergence of~${\cal W}$ when~$t$ tends to~$+\infty $ is thus sufficient to proof
the derivability of~${\cal X}(\theta ,p,p')$ with respect to~$p'$.

We have~${\cal Z}(\theta ,0,0,t)=1$ for all~$t,\theta $, and~${\cal Z}(\theta ,r,s,0)=1$ for all~$r,s,\theta $~;
thus~${\cal W}(\theta ,s,0)=0$ for all~$s,\theta $. The relation
\[
{\cal W}'(\theta ,s,t)=\frac{\partial}{\partial t}{\cal W}(\theta ,s,t)
=\frac{\partial}{\partial s} u_-(\psi_t^\kappa v(c(s),\theta )) 
\]
and the equation~\refeq{jacosta2} imply
\begin{equation}\label{eq:torsion2}
|{\cal W}'(\theta ,s,t)| \leq 
\Vert \nabla u_- \Vert_\infty C_1
 y_-(v(c(s),\theta ),t) .
\end{equation}
Thus we have for all~$t\geq 0$ and all~$s$~:
\[
|{\cal W}(\theta ,s,t)| \leq \Vert \nabla u_- \Vert_\infty C_1  
\int_{\tau =0}^t e^{-q_1\tau } \ d\tau 
\leq \frac{C_1\Vert \nabla u_- \Vert_\infty}{q_1}  ,
\]
\[
|{\cal Y}(\theta ,s,t)|\leq 
\frac{C_1\Vert \nabla u_- \Vert_\infty}{q_1}  y_-(v(c(s),\theta ),t),
\]
\[
\left| \frac{D}{ds} {\cal Z}(\theta ,0,s,t) \right| 
\leq \frac{C_1\Vert \nabla u_- \Vert_\infty}{q_1} {\cal Z}(\theta ,0,s,t)  .
\]
Integrating the last inequality with respect to~$s$ gives
\[
\exp \left( -\frac{C_1\Vert \nabla u_- \Vert_\infty}{q_1} | s | \right)  
\leq \left| {\cal Z}(\theta ,0,s,t) \right| \leq
\exp \left( \frac{C_1\Vert \nabla u_- \Vert_\infty}{q_1} | s | \right) .
\]
Because of~$| s |\leq L_0^s(c) \leq C_1 | s |$, we obtain
\[
\left| {\cal Z}(\theta ,0,s,t)-1 \right| 
\leq \exp \left( \frac{C_1\Vert \nabla u_- \Vert_\infty}{q_1}L_0^s(c)\right) 
-\exp \left( -\frac{C_1\Vert \nabla u_- \Vert_\infty}{q_1} L_0^s(c)\right) .
\]
Even by changing the horocycle, there comes for all~$t\geq 0$ and all~$r,s$~:
\[
\left| {\cal Z}(\theta ,r,s,t)-1 \right| 
\]
\[
\leq \exp \left( \frac{C_1\Vert \nabla u_- \Vert_\infty}{q_1}L_0^s(p_\theta (r,.))\right) 
-\exp \left( -\frac{C_1\Vert \nabla u_- \Vert_\infty}{q_1} L_0^s(p_\theta (r,.))\right) .
\]
We get for all~$t,\tau \geq 0$ and all~$s\in [-\varepsilon ,\varepsilon ]$~:
\[
\left| {\cal Z}(\theta ,0,s,t+\tau )-{\cal Z}(\theta ,0,s,t) \right| 
= \left| {\cal Z}(\theta ,0,s,t)\right| \left| {\cal Z}(\theta ,t,s,\tau )-1 \right| 
\]
\[
\leq \exp \left( \frac{C_1\Vert \nabla u_- \Vert_\infty}{q_1}L_0^s(c)\right) 
\]
\[
\times \left[ \exp \left( \frac{C_1\Vert \nabla u_- \Vert_\infty}{q_1}L_0^s(p_\theta (t,.))\right) 
-\exp \left( -\frac{C_1\Vert \nabla u_- \Vert_\infty}{q_1} L_0^s(p_\theta (t,.)\right) \right] .
\]
From the upper estimate~$L_0^s(p_\theta (r,.))\leq L_0^s(c)e^{-q_1t}$ results the existence of a constant~$C_3$
depending only from~$\Vert \nabla u_- \Vert_\infty$,~$\Vert \kappa  \Vert_\infty$,~$L_{-\varepsilon }^\varepsilon (c)$,~$q_1$ and~$q_0$
such that for all~$t,\tau $ positive and~$s\in [-\varepsilon ,\varepsilon ]$ we have
\begin{equation}\label{eq:fluconv}
\left| {\cal Z}(\theta ,0,s,t+\tau )-{\cal Z}(\theta ,0,s,t) \right| 
\leq C_3 L_0^s(p_\theta (t,.)) \leq C_3 L_{-\varepsilon }^\varepsilon (c) e^{-q_1t} .
\end{equation}
From the corollary~\ref{cor:courbure},
the equations~\refeq{decalage} and~\refeq{fluconv} imply~:
\[
\left| {\cal Z}(\theta ,r,s,t+\tau )-{\cal Z}(\theta ,r,s,t) \right| 
=\left| \frac{\displaystyle y_-(v(p_\theta (r,s),\theta ),t+\tau )}{\displaystyle y_-(v(p,\theta ),t+\tau )}
-\frac{\displaystyle y_-(v(p_\theta (r,s),\theta ),t)}{\displaystyle y_-(v(p,\theta ),t)}\right| 
\]
\[
\leq 
\frac{y_-(v(p_\theta (r,s),\theta ),-r)}{y_-(v(p_\theta (r,0),\theta ),-r)} C_3 f(d(c(-\varepsilon ),c(\varepsilon ))) e^{-q_1t} .
\]
The uniform Cauchy criterion on the compacts of~$\partial M\times M^2$
implies the convergence when~$t$ goes to~$+\infty $
of~${\cal Z}(\theta ,r,s,t)$ to a fonction~${\cal X}(\theta ,p,p_\theta (r,s))$ continuous on~$\partial M\times M^2$.

The equation~\refeq{decalage} ensures the existence and continuity of the partial derivative of~${\cal X}(\theta ,p,p_\theta (r,s))$ with respect to~$r$.

The equation~\refeq{torsion2} implies the upper estimate for~$t,\tau \geq 0$~:
\[
\left| {\cal W}(\theta ,s,t+\tau ) - {\cal W}(\theta ,s,t ) \right| 
\leq \int_{\rho =t}^{t+\tau } C_1 e^{-q_1\rho }  
\left\Vert \nabla u_- \right\Vert_\infty \ d\rho 
\]
thus
\begin{equation}\label{eq:derifluc}
\left| {\cal W}(\theta ,s,t+\tau ) - {\cal W}(\theta ,s,t ) \right| 
\leq e^{-q_1t}\frac{C_1}{q_1}  \left\Vert \nabla u_- \right\Vert_\infty ,
\end{equation}
which ensures the convergence, uniform over the compacts of~$\partial M\times \R $, of the family of mappings~${\cal W}(\theta ,s,t)$,
when~$t$ goes to~$+\infty $, to a continuous mapping of~$\theta ,s$.
The relations~\refeq{decalage},~\refeq{derifluc} and
\[ 
\frac{\partial }{\partial s} {\cal Z}(\theta ,0,s,t) = {\cal W}(\theta ,s,t){\cal Z}(\theta ,0,s,t)
\]
imply the uniform convergence, over the compacts of~$\partial M\times \R ^2$, 
of the family of mappings~$\partial {\cal Z}/\partial s (\theta ,r,s,t) $,
when~$t$ goes to~$+\infty $, to a continuous mapping.
This gives the condition of derivability with respect to~$s$ for the function~${\cal X}(\theta ,p,p_\theta (r,s))$,
thus the derivability with respect to~$p'$ anounced for the function~${\cal X}(\theta ,p,p')$.
The trivial relation
\begin{equation}\label{flucinv}
{\cal X}(\theta ,p,p') \cdot {\cal X}(\theta ,p',p)=1
\end{equation}
gives the derivability with respect to~$p$.
\end{preuve}

\begin{definition}\label{def:transfert}
With the assumptions of the theorem~\ref{thm:fluc}, for~$v\in T^1M$,~$v'\in W^{CS}(v)$, 
the limit mapping calculated in the theorem~\ref{thm:fluc} is called {\em stable transfer} from ~$v'$ to~$v$ and is written down~:
\[
X(v,v')=\lim_{t\to +\infty } \frac{y_-(v',t)}{y_-(v,t)} = {\cal X}(v_{+\infty }, \pi v,\pi v').
\]
The~{\em extended stable transfer} from~$v'$ to~$v$ is the mapping which to~$\xi \in T_{\pi v'}M$ associates
\[
\widetilde{X} (v,v')\xi =X(v,v') \langle \xi ,N(v')\rangle N(v) + \langle \xi ,v'\rangle v.
\]
\end{definition}

\begin{cor}\label{cor:regutransfert}
Let~$M$ be a complete, connected, oriented, simply connected surface, equipped with a~${\cal C}^\infty $-riemannian metric~$g$
whose curvature is negatively pinched~:~$-k_0^2\leq K\leq -k_1^2<0$, and with a~${\cal C}^\infty $-magnetic field~$\kappa $
whose~${\cal C}^1$-norm is bounded, with Jacobi endomorphism satisfying the pinching condition~\refeq{pince},
and such that the gradient of the centre-stable space~$u_-$ is uniformly bounded.
The stable transfer and extended stable transfer are of~${\cal C}^1$-class on a given centre-stable manifold.
\end{cor}
\begin{preuve}
This follows from the derivability of the stable transfer stated in the theorem~\ref{thm:fluc}.
\end{preuve}

The symplectic structure of the space of geodesics (section~\ref{sec:symp}) leads to the following result.

\begin{thm}\label{thm:croisefluc}
Let~$M$ be a complete, connected, oriented, simply connected surface, equipped with a~${\cal C}^\infty $-riemannian metric~$g$
whose curvature is negatively pinched~:~$-k_0^2\leq K\leq -k_1^2<0$, and with a~${\cal C}^\infty $-magnetic field~$\kappa $
whose~${\cal C}^1$-norm is bounded, with Jacobi endomorphism satisfying the pinching condition~\refeq{pince},
and such that the gradient of the centre-stable~$u_-$ (respectively centre-unstable~$u_+$) spaces is uniformly bounded.
For a given point~$\theta $ at infinity, and  unitary vectors~$v$,~$v'$ belonging to the centre-stable manifold~$W^{CS}(\theta )$ 
determined by~$\theta $, the family of mappings~$(v,v') \mapsto y_+(v',t)/y_+(v,t)$
converges uniformly when~$t$ goes to~$+\infty $ to a continuous mapping.
The limit mapping admits a partial derivative with respect to~$v'$ (in the direction~$W^{CS}(\theta )$),
continuous with respect to~$(v,v')$.
\end{thm}
\begin{preuve}
Following the formula~\refeq{sympl}, for~$v,v'\in T^1M$,~$t\in \R^+$, we have
\[
\frac{y_+(v',t)}{y_+(v,t)} = 
\frac{(u_+-u_-)(v')\cdot (u_+-u_-)(\psi ^\kappa _tv)\cdot y_-(v,t)}{(u_+-u_-)(\psi ^\kappa _tv')\cdot y_-(v',t)\cdot (u_+-u_-)(v)} .
\]
Because the gradient of~$(u_+-u_-)$ is uniformly bounded, the quotient~$(u_+-u_-)(\psi ^\kappa _tv)/(u_+-u_-)(\psi ^\kappa _tv')$
tends to one uniformly over the compacts, when ~$t$ goes to~$+\infty $, this fact implying the uniform convergence over the compacts~:
\[
\frac{y_+(v',t)}{y_+(v,t)} \quad \mathop{\longrightarrow}_{t\to+\infty } \quad \frac{(u_+-u_-)(v')}{(u_+-u_-)(v)} X(v',v) .
\]
The regularity of the limit results from the theorem~\ref{thm:fluc}.
\end{preuve}

\begin{definition}\label{def:transfU}
With the assumptions of the theorem~\ref{thm:croisefluc}, for~$v\in T^1M$,~$v'\in W^{CS}(v)$, 
the limit mapping calculated in the theorem~\ref{thm:croisefluc} is called {\em unstable transfer}
from~$v'$ to~$v$ and is written down~:
\[
\underline{X} (v,v') =\lim_{t\to +\infty } \frac{y_+(v',t)}{y_+(v,t)} 
=\frac{(u_+-u_-)(v')}{(u_+-u_-)(v)} X(v',v) .
\]
The~{\em extended unstable transfer} from~$v'$ to~$v$ is the mapping which to~$\xi \in T_{\pi v'}M$ associates
\[
\widetilde{\underline{X}} (v,v')\xi =\underline{X}(v,v') \langle \xi ,N(v')\rangle N(v) + \langle \xi ,v'\rangle v.
\]
\end{definition}

The following result is immediate.

\begin{cor}\label{cor:regutransU}
Let~$M$ be a complete, connected, oriented, simply connected surface, equipped with a~${\cal C}^\infty $-riemannian metric~$g$
whose curvature is negatively pinched~:~$-k_0^2\leq K\leq -k_1^2<0$, and with a~${\cal C}^\infty $-magnetic field~$\kappa $
whose~${\cal C}^1$-norm is bounded, with Jacobi endomorphism satisfying the pinching condition~\refeq{pince},
and such that the gradient of the centre-stable~$u_-$ (respectively centre-unstable~$u_+$) spaces is uniformly bounded.
The unstable transfer and the extended unstable transfer are of ~${\cal C}^1$-class on a given centre-stable manifold.
\end{cor}

%%%%%%%%%%%%%%%%%%%%%%%%%%%%%%%%%%%%%%%%%%%%%%%%%%%

\section{Lyapounoff exponents}\label{sec:lyap}

%%%%%%%%%%%%%%%%%%%%%%%%%%%%%%%%%%%%%%%%%%%%%%%%%%%

\begin{thm}\label{thm:Lyap}
Let~$M$ be a compact, connected, oriented, surface, equipped with a~${\cal C}^\infty $-riemannian metric~$g$
whose curvature is negatively pinched~:~$-k_0^2\leq K\leq -k_1^2<0$, and with a~${\cal C}^\infty $-magnetic field~$\kappa $
with Jacobi endomorphism satisfying the pinching condition~\refeq{pince}.
If the (future) Lyapounoff exponents of the magnetic flow are defined at an element~$v$ of~$T^1M$, 
then they are defined and constant along its centre-stable manifold~$W^{CS}(v)$.
\end{thm}

\begin{preuve}
The compacity insures the uniform bounds over the gradients of~$\kappa $,~$u_+$ and~$u_-$. 
It is sufficient to pass to the universal cover of~$M$ and to apply the theorems~\ref{thm:fluc} and~\ref{thm:croisefluc}.
\end{preuve}

%%%%%%%%%%%%%%%%%%%%%%%%%%%%%%%%%%%%%%%%%%%%%%%%%%%

\section{Horocyclic transport}\label{sec:transport}

%%%%%%%%%%%%%%%%%%%%%%%%%%%%%%%%%%%%%%%%%%%%%%%%%%%

In this section are collected some tools for the section~\ref{sec:linea}~;
the notations and the assumptions are those of the section~\ref{sec:fluc},
particularly of the proof of the theorem~\ref{thm:fluc}.
Let~$\tau (\theta ,s,t)$ be the parallel transport along the curve~$p_\theta (t,.)$ which sends~$T_{p_\theta (t,s)}M$ onto~$T_{p_\theta (t,0)}M$. 
Let~$\zeta (\theta ,s,t)$ be the angle between the vectors~$\tau (\theta ,s,t)\cdot T_\theta (s,t)$ and~$T_\theta (0,t)$. 
Composing~$\tau $ with the rotation of angle~$\zeta $ yields the isometry
\[
\chi (\theta ,s,t)=\rho _{\zeta (\theta ,s,t)} \tau (\theta ,s,t)
\]
which sends the direct orthonormal basis~$(T_\theta (s,t),N_\theta (s,t))$ onto the direct orthonormal basis~$(T_\theta (0,t),N_\theta (0,t))$.
\begin{definition}
The mapping~$\chi (\theta ,s,t)$ is called {\em horocyclic transport}.
\end{definition}
\begin{rem}
The horocyclic transport is continuous on~$\partial M\times\R^2$.
\end{rem}
In this section the control of~$\chi $ is precised in different ways.

\begin{lem}\label{lem:transport1}
With the above notations, we have for all~$\theta \in \partial M$, $s\in [-\varepsilon ,\varepsilon ]$, $t\in \R^+$~:
\[
\left\Vert \frac{D}{dt} \chi (\theta ,s,t) \right\Vert 
\leq \left( 2k_0^2 + \left\Vert \nabla \kappa \right\Vert_\infty \right) C_1 e^{-q_1t} | s |.
\]
\end{lem}

\begin{preuve}
For all~$\theta \in \partial M$,~$t\in \R$,~$s\in]-\varepsilon ,\varepsilon [$,~$\xi \in T_{p_\theta (t,s)}M$,~$\eta \in T_{p_\theta (t,0)}M$, we have
\[
\left\{
\matrix{
\displaystyle \frac{D}{ds} \tau (\theta ,s,t)\xi =0 , \cr 
\tau (\theta ,0,t)\eta =\eta .
}
\right.
\]
The field~$D\tau (\theta ,s,t)/dt$ is the solution~${\cal T}$ of
\[
\left\{
\matrix{
\displaystyle \frac{D{\cal T}}{ds} (\theta ,s,t) = R(T,S) \tau (\theta ,s,t) , \cr 
{\cal T} (\theta ,0,t) = 0 ,
}
\right.
\]
thus, due to the upper estimate~\refeq{jacosta2}, it satisfies the differential inequality
\[
\left\Vert \frac{D{\cal T}}{ds} (\theta ,s,t) \right\Vert \leq 
k_0^2 C_1 e^{-q_1t}  ,
\]
therefore
\begin{equation}\label{eq:trans0}
\left\Vert \frac{D\tau }{dt} (\theta ,s,t) \right\Vert \leq 
\int_{\sigma =0}^s k_0^2 C_1 e^{-q_1t} \ d\sigma 
\leq k_0^2 C_1 e^{-q_1t} | s |.
\end{equation}
We have
\[
\cos \zeta (\theta ,s,t) =\left\langle \tau (\theta ,s,t) T_\theta (s,t) , T_\theta (0,t) \right\rangle ,
\]
\[
\sin \zeta (\theta ,s,t) =\left\langle \tau (\theta ,s,t) T_\theta (s,t) , N_\theta (0,t) \right\rangle ,
\]
thus
\[
\frac{\partial }{\partial t} \cos \zeta (\theta ,s,t) 
= \left\langle \frac{D\tau }{dt} (\theta ,s,t) T_\theta (s,t) , T_\theta (0,t) \right\rangle 
\]
\[
+ \left\langle \tau (\theta ,s,t) \kappa (p_\theta (t,s)) N_\theta (s,t) , T_\theta (0,t) \right\rangle 
+ \left\langle \tau (\theta ,s,t) T_\theta (s,t) , \kappa (p_\theta (t,0)) N_\theta (0,t) \right\rangle 
\]
\[
= \left\langle \frac{D\tau }{dt} (\theta ,s,t) T_\theta (s,t) , T_\theta (0,t) \right\rangle 
+\sin \zeta (\theta ,s,t) \left[ \kappa (p_\theta (t,0)) - \kappa (p_\theta (t,s)) \right] ,
\]
and
\[
\frac{\partial }{\partial t} \sin \zeta (\theta ,s,t) 
= \left\langle \frac{D\tau }{dt} (\theta ,s,t) T_\theta (s,t) , N_\theta (s,t) \right\rangle 
\]
\[
+ \left\langle \tau (\theta ,s,t) \kappa (p_\theta (t,s)) N_\theta (s,t) , N_\theta (0,t) \right\rangle 
+ \left\langle \tau (\theta ,s,t) T_\theta (s,t) , -\kappa (p_\theta (t,0)) T_\theta (0,t) \right\rangle 
\]
\[
= \left\langle \frac{D\tau }{dt} (\theta ,s,t) T_\theta (s,t) , N_\theta (0,t) \right\rangle 
+\cos \zeta (\theta ,s,t) \left[ \kappa (p_\theta (t,s)) - \kappa (p_\theta (t,0)) \right] .
\]
We deduce from this
\[
\frac{\partial \zeta }{\partial t} (\theta ,s,t) =\kappa (p_\theta (t,s)) - \kappa (p_\theta (t,0)) 
- \left\langle \frac{D\tau }{dt} (\theta ,s,t) T_\theta (s,t) , T_\theta (0,t) \right\rangle 
\sin \zeta (\theta ,s,t) 
\]
\[
+ \left\langle \frac{D\tau }{dt} (\theta ,s,t) T_\theta (s,t) , N_\theta (0,t) \right\rangle 
\cos \zeta (\theta ,s,t) .
\]
The formula~\refeq{jacosta2} implies~:
\[
\left| \kappa (p_\theta (t,s)) - \kappa (p_\theta (t,0)) \right| \leq \left\Vert \nabla \kappa \right\Vert_\infty C_1 e^{-q_1t} | s | .
\]
Following, due to the upper estimate~\refeq{trans0}, the rotation satisfies the differential inequality
\[
\left\Vert \frac{D}{dt} \rho _{\zeta (\theta ,s,t)} \right\Vert 
\leq \left| \frac{\partial}{\partial t} \zeta (\theta ,s,t) \right| 
\leq \left\Vert \frac{D\tau }{dt} (\theta ,s,t) \right\Vert + \left| \kappa (p_\theta (t,s)) - \kappa (p_\theta (t,0)) \right| 
\]
\[
\leq \left( k_0^2 + \left\Vert \nabla \kappa \right\Vert_\infty \right) C_1 e^{-q_1t} | s |.
\]
The definition of~$\chi $ permits to conclude.
\end{preuve}

We have also the upper estimate
\[
\left\Vert \frac{D}{ds} \chi (\theta ,s,t) \right\Vert \leq 
\left\Vert \frac{D}{ds} \rho _{\zeta (\theta ,s,t)} \right\Vert \leq 
\left\Vert \frac{D}{ds} T_\theta (s,t) \right\Vert =
\left\Vert \frac{D}{dt} S_\theta (s,t) \right\Vert 
\]
\[
\leq C_1 e^{-q_1t}  ,
\]
from which follows the~:
\begin{lem}\label{lem:transport2}
With the preceding notations, we have for all~$\theta \in \partial M$, $s\in [-\varepsilon ,\varepsilon ]$, $t\in \R^+$~:
\[
\left\Vert \frac{D}{ds} \chi (\theta ,s,t) \right\Vert 
\leq C_1 e^{-q_1t} \leq C_1 e^{-q_1t} \left\Vert \frac{dc}{ds}(s) \right\Vert .
\]
\end{lem}
The horocyclic transport presents some uniformity.
\begin{lem}\label{lem:transport3}
Let~$A:v\in T^1M\mapsto A(v)\in L(T_{\pi (v)}M)$ be a field of linear endomorphisms of class~${\cal C}^1$ over~$T^1M$, bounded
in~${\cal C}^1$-norm. With the preceding notations,
for all~$\theta \in \partial M$,~$s\in [-\varepsilon ,\varepsilon ]$ and~$t\in\R ^+$, we have
\[
\left\Vert A(v(p_\theta (0,t),\theta ))\chi (\theta ,s,t)-\chi (\theta ,s,t) A(v(p_\theta (s,t),\theta )) \right\Vert 
\]
\[
\leq  C_1 \left( 2 \left\Vert A \right\Vert _\infty +\left\Vert DA \right\Vert _\infty \right) e^{-q_1t} | s | .
\] 
\end{lem}
\begin{preuve}
We have
\[
\left\Vert \frac{D \chi (\theta ,s,t) A(v(p_\theta (s,t),\theta ))}{ds} \right\Vert
\leq \left\Vert A(v(p_\theta (s,t),\theta )) \right\Vert \left\Vert \frac{D}{ds} \chi (\theta ,s,t) \right\Vert
\]
\[+ \left\Vert \chi (\theta ,s,t) \right\Vert \left\Vert \frac{D}{ds} A(v(p_\theta (s,t),\theta )) \right\Vert ,
\]
which is bounded from above, following the lemma~\ref{lem:transport2}, by
\[
\left\Vert A \right\Vert _\infty C_1 e^{-q_1t}  
+ \left\Vert DA \right\Vert _\infty 
\left( \left\Vert S(s,t) \right\Vert +\left\Vert \frac{D}{dt} S(s,t) \right\Vert \right) 
\]
\[
\leq \left( \left\Vert A \right\Vert _\infty +\left\Vert DA \right\Vert _\infty \right) 
C_1 e^{-q_1t}  .
\]
The covariant derivative with respect to~$s$ of~$A(v(p_\theta (0,t),\theta )) \chi (\theta ,s,t) $
admits the same first term as above for upper bound.
The quantity that we aim to estimate in the lemma annulates at~$s=0$,
the upper estimate of the statement is obtained by integrating the expression
\[
 \left( 2 \left\Vert A \right\Vert _\infty +\left\Vert DA \right\Vert _\infty \right) 
C_1 e^{-q_1t}  
\]
over~$s$ from~$0$.
\end{preuve}

%%%%%%%%%%%%%%%%%%%%%%%%%%%%%%%%%%%%%%%%%%%%%%%%%%%

\section{Linearization}\label{sec:linea}

%%%%%%%%%%%%%%%%%%%%%%%%%%%%%%%%%%%%%%%%%%%%%%%%%%%

\begin{definition}\label{def:poussee}
Let~$M$ be a complete, connected, oriented, simply connected surface, equipped with a~${\cal C}^\infty $-riemannian metric~$g$
whose curvature is negatively pinched~:~$-k_0^2\leq K\leq -k_1^2<0$, and with a~${\cal C}^\infty $-magnetic field~$\kappa $
whose~${\cal C}^1$-norm is bounded, with Jacobi endomorphism satisfying the pinching condition~\refeq{pince},
and such that the gradient of the centre-stable space~$u_-$ is uniformly bounded.
For~$v\in T^1M$,~$t\in \R$, the {\em stable push} is the mapping~$\Phi ^{\kappa ,v}_t:M\to M$ which to every~$p\in M$ 
associates~$\Phi ^{\kappa ,v}_tp=\pi \psi ^\kappa _tv^\kappa (p,v^\kappa _{+\infty})$.
\end{definition}

\begin{thm}\label{thm:linea}
With the assumptions of the definition~\ref{def:poussee},
for all unitary vector~$v\in T^1M$, there exists a unique mapping~$E_v$ 
from~$M$ into~$T_{\pi (v)}M$ such that~:
\[
\begin{matrix}
{
i) & \forall z\in M \quad \left\langle E_v(z),v\right\rangle = B_v(z) , \cr
{} & {} \cr
ii) & \forall z\in M , \forall t\in \R \quad E_v\left( \Phi^{\kappa ,v}_t (z)\right) = tv+E_v(z) , \cr
{} & {} \cr
iii) & \forall z \in H^\kappa _v(0) \quad E_v(z) = 
{\displaystyle \lim_{t \to +\infty } y_-(v,t) ^{-1} 
\left\langle \left( \exp _{\pi \psi ^\kappa _t v} \right) ^{-1} \Phi _t^{\kappa ,v}(z) , n \right\rangle n,
 } \cr 
}
\end{matrix}
\]
writing down~$n= N(\psi ^\kappa _tv)$. Moreover, the mapping~$(v,p)\mapsto E_v(p)$ is continuous from~$T^1M\times M$ into~$TM$.
\end{thm} 

\begin{preuve}
The conditions~$i)$,~$ii)$ and~$iii)$ ensure naturally the uniqueness of the solution.
Let~$v$ be in~$T^1M$. Let us define~$\theta =v_{+\infty }$,~$p=\pi v$.
Taking the notations of the sections~\ref{sec:fluc} and~\ref{sec:transport}, 
we have~$N_\theta (0,t)=n$ and we define~:
\[
X_{\theta ,s} = \int_{\sigma =0}^s \chi (\theta ,\sigma ,0)\frac{dc}{ds}(\sigma ) \ d\sigma ,
\qquad X_{\theta ,s}^t = \int_{\sigma =0}^s \chi (\theta ,\sigma ,t)\frac{\partial p_\theta (t,\sigma )}{\partial s}(\sigma ) \ d\sigma ,
\]
\[
e_v^t(s)=\frac{\left\langle n , X_{\theta ,s}^t \right\rangle }{y_-(v,t)} .
\]
We have
\[
\chi (\theta ,s,t) T_\theta (s,t)=\psi ^\kappa _tv, \qquad \chi (\theta ,s,t) N_\theta (s,t)=N(\psi ^\kappa _tv)=N_\theta (0,t) ,
\]
thus
\[
e_v^t(s)=\frac{\displaystyle \int_{\sigma =0}^s 
\hskip -10 pt
y_-(v(p_\theta (t,\sigma ),\theta ) ,t)\left\langle \frac{dc}{ds}(\sigma ) , N_\theta (\sigma ,0) \right\rangle \ d\sigma }{y_-(v,t)} 
\]
\[
=\int_{\sigma =0}^s 
\hskip -10 pt
{\cal Z}(\theta ,0,\sigma ,t)\left\langle \frac{dc}{ds}(\sigma ) , N_\theta (\sigma ,0) \right\rangle \ d\sigma .
\]
The equation~\refeq{fluconv} implies for~$t,\tau \geq 0$:
\[
\left| e_v^{t+\tau }(s)-e_v^t (s) \right| \leq 
\hskip -3 pt
\int_{\sigma \in [0,s]} 
\hskip -20 pt
C_3 L_{-\varepsilon }^\varepsilon (c) e^{-q_1t}
\left| \left\langle \frac{dc}{ds}(\sigma ) , N_\theta (\sigma ,0) \right\rangle \right|d\sigma 
\leq C_3 L_{-\varepsilon }^\varepsilon (c)^2 e^{-q_1t} .
\]
By taking 
\[
{\cal E}^t_v(p_\theta (r,s))=rv+e^t_v(s)N(v) ,
\]
and with the foregoing upper estimate and the corollary~\ref{cor:courbure} we get~:
\[
\left|{\cal E} _v^{t+\tau }(p_\theta (r,s))-{\cal E}_v^t (p_\theta (r,s)) \right| \leq 
C_3 L_{-\varepsilon }^\varepsilon (c)^2 e^{-q_1t} \leq C_3 f\left( d(c(-\varepsilon ),c(\varepsilon ) )\right) ^2 e^{-q_1t} .
\]
For every compact~${\cal K}$ of~$T^1M\times M$, the family of continuous mappings~${\cal E}_v^t(p_\theta (r,s))$
satisfies the uniform Cauchy condition over~${\cal K}$ when~$t$ tends to~$+\infty $.
It converges thus towards a mapping continuous over~$T^1M\times M$.
Establishing the formulas~$i)$ and~$ii)$ is immediate.
There remains to proof the formula~$iii)$.

For~$w\in T^1M$,~$r\in \R$,~$\xi \in T_{\pi w}M$, if~$\tilde J(w,r)$ is the {\em geodesic} Jacobi field
along the geodesic curve directed by~$w$, such that~$\tilde J(w,0)=0_{T_{\pi w}M}$
and~$\tilde J'(w,0)=\id_{T_{\pi w}M}$, the linear mapping tangent to the exponential satisfies~(\cite{Ga-Hu-La},~3.46 p.117)~:
\begin{equation}\label{eq:bidule2} 
T_\xi \exp_{\pi w} rw = \frac1r \tilde J(w,r) \xi .
\end{equation}
From the bounds on the Gauss cuvature results the existence of a constant~$C_4$ such that for all~$r\in [0,1]$~:
\[
\left\Vert \tilde J(w,r) -r \id_{T_{\pi w}M} \right\Vert \leq C_4 r^3 .
\]
The derivative of the exponential in every zero vector is the identity of
the tangent  vector space, and~$\chi (\theta ,0,t)$ is the identity of~$T_{p_\theta (t,0)}M$. 
Even by situating in a chart in the neighbourhood of~$p_\theta (t,s)$, 
there exists a constant~$C_5$ such that~:
\begin{equation}\label{eq:bidule3}
\left| \left\langle \exp_{\pi \psi ^\kappa _tv}^{-1} p_\theta (t,s) ,n \right\rangle -y_-(v,t)e_v^t(s) \right|
\]
\[
\leq \int_{\sigma =0}^s
\left\Vert \frac{\partial }{\partial s} \exp_{\pi \psi ^\kappa _tv}^{-1} p_\theta (t,\sigma ) 
- \chi (\theta ,\sigma ,t)\frac{\partial p_\theta  }{\partial s} (t,\sigma ) \ d\sigma  \right\Vert 
\]
\[
\leq 
C_5 d\left( \pi \psi ^\kappa _tv, p_\theta (t,s) \right) ^2 
+C_1 e^{-q_1t} \int_{\sigma '=0}^s\int_{\sigma =0}^{\sigma '} 
\left\Vert \frac{\partial p_\theta }{\partial s} (t,\sigma ) \right\Vert 
\left\Vert \frac{dc }{ds} (\sigma ') \right\Vert \ d\sigma d\sigma '
\]
\[
\leq 
C_5 d\left( \pi \psi ^\kappa _tv, p_\theta (t,s) \right) ^2 
+ C_5 L_0^\varepsilon (p_\theta (t,.)) e^{-q_1t} L_{-\varepsilon }^\varepsilon (c).
\end{equation}
The first term of the last member comes from the effect of the bounds over the Gauss curvature on the exponential, 
and the second term comes from the lemma~\ref{lem:transport2} by carrying out two successive integrations.
We deduce the inequality~:
\[
\left| \frac{ \displaystyle \left\langle \exp_{\pi \psi ^\kappa _tv}^{-1} p_\theta (t,s) ,n \right\rangle }{y_-(v,t)} -e_v^t(s) \right|
\]
\[
\leq C_5 \left[ d\left( \pi \psi ^\kappa _tv, p_\theta (t,s) \right) + e^{-q_1t} L_{-\varepsilon }^\varepsilon (c) \right]
\int_{\sigma \in [0,s]} 
\hskip -20 pt
\frac{y_-(v(p_\theta (t,\sigma ),\theta ),t)}{y_-(v,t)} 
\left\Vert \frac{dc}{ds}(\sigma ) \right\Vert \ d\sigma .
\]
The quantity~$d\left( \pi \psi ^\kappa _tv, p_\theta (t,s) \right)$ is bounded from above 
by~$L_0^s(p_\theta (t,.)$ which is inferior or equal to~$C_1 L_{-\varepsilon }^\varepsilon (c) e^{-q_1t} $.
According to the theorem~\ref{thm:fluc}, there exists a constant~$C_6$ (depending on~$\varepsilon $)
such that the last integral above is bounded by~$C_6 L_0^s(c)$ independently of~$t$,
from which we get the upper estimate~:

\[
\left| \frac{ \displaystyle \left\langle \exp_{\pi \psi ^\kappa _tv}^{-1}  p_\theta (t,s) ,n \right\rangle }{y_-(v,t)} -e_v^t(s) \right|
\leq (1+C_1) C_5 C_6 L_{-\varepsilon }^\varepsilon (c)^2 e^{-q_1t} .
\]
The limit stated in the formula~$iii)$ of the theorem is thus obtained.
\end{preuve}

The linearization cooperates to some extent with the magnetic flow.

\begin{cor}\label{cor:linea}
For~$v\in T^1M$,~$p\in M$,~$t\in\R$, we have
\[
E_{\psi _tv}(p)=\left( \big\langle E_v(p),v\big\rangle -t \right) \psi _t v
+ y_-(v,t)\big\langle E_v(p),N(v) \big\rangle N(\psi _tv) .
\]
\end{cor}
\begin{preuve}
This immediately results from the construction of the theorem~\ref{thm:linea}.
\end{preuve}

%%%%%%%%%%%%%%%%%%%%%%%%%%%%%%%%%%%%%%%%%%%%%%%%%%%

\section{Regularity of the linearization}\label{sec:regu}

%%%%%%%%%%%%%%%%%%%%%%%%%%%%%%%%%%%%%%%%%%%%%%%%%%%

\begin{thm}\label{thm:regulin}
Let~$M$ be a complete, connected, oriented, simply connected surface, equipped with a~${\cal C}^\infty $-riemannian metric~$g$
whose curvature is negatively pinched~:~$-k_0^2\leq K\leq -k_1^2<0$, and with a~${\cal C}^\infty $-magnetic field~$\kappa $
whose~${\cal C}^1$-norm is bounded, with Jacobi endomorphism satisfying the pinching condition~\refeq{pince},
and such that the gradient of the centre-stable~$u_-$ 
%%
%%(respectively centre-unstable~$u_+$) 
%%
spaces is uniformly bounded.
For~$v\in T^1M$ and~$\theta =v^\kappa _{+\infty }$, the mapping~$E_v$ is of~${\cal C}^2$ over~$M$, 
and with the extended stable transfer~$\widetilde{X}$ coming from the definition~\ref{def:transfert},
its derivative at~$y\in M$ is~:
\[
\widetilde{X}\left( v,v^\kappa (y,\theta ) \right) .
\]
\end{thm}

\begin{preuve}
With the notations of the sections~\ref{sec:fluc} and~\ref{sec:linea}, we have
\[
\frac{\partial e_v^t(s)}{\partial s}={\cal Z}(\theta ,0,s,t) ,
\]
which converges when~$t$ tends to~$+\infty $ towards
\[
X(v,v_s) .
\]
The longitudinal component of the derivative of~$E_v$ is~$v$, and the continuous derivability of the derivative of~$E_v$ results from the corollary~\ref{cor:regutransfert}.
\end{preuve}

\begin{cor}\label{cor:diffeo}
With the same assumptions as in the theorem~\ref{thm:regulin},
for all~$v\in T^1M$, the mapping~$E_v$ is a~${\cal C}^2$-diffeomorphism.
\end{cor}

\begin{preuve}
It is clear that the mapping~$E_v$ is surjective.
It is a local~${\cal C}^2$-diffeomorphism according to the theorem~\ref{thm:regulin}, 
it is therefore a cover of~$T_{\pi (v)}M$ (which is isomorphic to~$\R ^2$), 
and consequently it is a diffeomorphism.
\end{preuve}

\begin{rem}
We may notice that the mapping~$v\mapsto E_v(p)$ with~$p$ fixed cannot be absolutely continuous.
\end{rem}
Let~$v$ and~$w$ be two unitary tangent vectors belonging to a same unstable manifold,
such that~$p$ belongs to the curve directed by~$v$.
The vectors ~$\psi_tv$ and~$\psi_tw$ are arbitrarily close when~$t$ tends to~$-\infty $.
The orthogonal component of~$E_{\psi_tv}(p)$ is zero and the norm of the orthogonal component of~$E_{\psi_tw}(p)$ 
tends to~$+\infty $ when~$t$ tends to~$-\infty $.

%%%%%%%%%%%%%%%%%%%%%%%%%%%%%%%%%%%%%%%%%%%%%%%%%%%

\section{Flow conjugacy}\label{sec:conju}

%%%%%%%%%%%%%%%%%%%%%%%%%%%%%%%%%%%%%%%%%%%%%%%%%%%

\begin{thm}\label{thm:reguconj}
Let~$M$ be a closed (compact without boundary), connected, oriented surface.
Let~$g_1$ and~$g_2$ be two~${\cal C}^\infty $-riemannian metrics over~$M$
whose curvatures are negatively pinched~:~$-k_0^2\leq K_i\leq -k_1^2<0$ for~$i=1,2$.
Let~$\kappa _1$,~$\kappa _2$ be two~${\cal C}^\infty $-magnetic fields over~$M$.
The two magnetic flows~$\psi ^1_t=\psi ^{g_1,\kappa _1}_t$ and~$\psi ^2_t=\psi ^{g_2,\kappa_2}_t$ 
are supposed to have negative Jacobi endomorphisms (thus satisfying the pinching condition~\refeq{pince}).
If the two magnetic flows have the same marked length spectrum,
then they are conjugated by a bi-h\"older-continuous homeomorphism~$h$ from~$T^1_1M$ onto~$T^1_2M$ 
which preserves the Lyapounoff exponents of periodic orbits.
\end{thm}

\begin{preuve}
The gradient of the centre-stables spaces~$u_-$ and of the centre-instables spaces~$u_+$ 
are uniformly bounded because the unitary tangent bundles are compact.
The Jacobi endomorphisms satisfie the pinching condition~\refeq{pince} for the same reason.
The existence of the bi-h\"older-continuous conjugacy~$h$ from~$T^1_1M$ onto~$T^1_2M$ is well-known~\cite{EHS,Ha2,Li}. 
Let~$v$ be a~$T$-periodic vector for~$\psi ^1$~; its conjugate ~$hv$ is~$T$-periodic for~$\psi ^2$
and the conjugacy~$h$ maps~$W^S(v)$ onto~$W^S(hv)$. The linearizations associated to~$\psi ^i$ are written down~$E^i$ for~$i=1,2$. 
Let~${\cal A}$ be the bijection, restricted to the orthogonal spaces identified to the real line, defined as follows~:
\[
\matrix{
{\cal A} : & \R \simeq \R N_1(v)  & \hskip -3 pt \displaystyle \mathop{\longrightarrow}^{\displaystyle \left( E^1_v \circ \pi \right) ^{-1}} & \hskip -9 pt W^S(v) 
& \hskip -3 pt \displaystyle \mathop{\longrightarrow}^{\displaystyle h} & \hskip -3 pt W^S(hv) 
& \hskip -3 pt \displaystyle \mathop{\longrightarrow}^{\displaystyle E^2_{hv} \circ \pi } & \hskip -3 pt \R N_2(hv)\simeq \R     \cr
{}         & \xi N_1(v)           & \hskip -3 pt \longmapsto                                                                                 & \hskip -9 pt w
& \hskip -3 pt \longmapsto                                              & \hskip -3 pt hw      
& \hskip -3 pt \longmapsto                                                                & \hskip -3 pt {\cal A}(\xi )N_2(hv) . \cr
}
\]
The Lyapounoff exponents are written down~:
\[
\lambda _{-,1}(v)=-\lambda _{+,1}(v) \quad \hbox{\rm and} \quad \lambda _{-,2}(hv)=-\lambda _{+,2}(hv).
\]
Let us denote~:
\[
\nu _1=e^{T\lambda _{-,1}(v)}=y_{-,1}(v,T)\quad \hbox{\rm and} \quad \nu _2=e^{T\lambda _{-,2}(hv)}=y_{-,2}(hv,T).
\]
According to the corollary~\ref{cor:linea}, for each of the two flows, every~$T$-periodic vector~$v$ satisfies~:
\[
E_{\psi_{{}_T}v}\Big| _{\displaystyle W^S(v)}=y_-(v,T)E_v\Big| _{\displaystyle W^S(v)} .
\]
For all~$n\in\N$ and~$\xi \in \R$ we therefore have~:
\[
{\cal A}(\nu _1^n\xi )=\nu _2^n{\cal A}(\xi ) .
\]
Because the conjugacy~$h$ is bi-h\"older-continuous and because the linearizations are of~${\cal C}^2$-class,
there exist two constants~$C>0$ and~$\alpha \in ]0,1]$ such that for all~$\xi \in [-1,1]$ we have~:
\[
\left| {\cal A}(\xi ) \right| \leq C |\xi |^\alpha \quad \hbox{\rm and} \quad |\xi |\leq C \left| {\cal A}(\xi ) \right| ^\alpha .
\]
Thus for all natural integer~$n$ we have~:
\[
\nu _2^n |{\cal A}(1)| \leq C \nu _1^{n\alpha }\quad \hbox{\rm and} \quad \nu _1^n\leq C |{\cal A}(1)|^\alpha \nu _2^{n\alpha } ,
\]
which implies for all~$n\in \N ^\star $~:
\[
\ln \nu _2 \leq \alpha \ln \nu _1 +\frac{\ln C-\ln |{\cal A}(1)|}{n} \quad \hbox{\rm and} \quad 
\ln \nu _1 \leq \alpha \ln \nu _2 +\frac{\ln C+\alpha \ln |{\cal A}(1)|}{n} .
\]
By making~$n$ tend to~$+\infty $ there comes
\[
\ln \nu _2 \leq \alpha \ln \nu _1 \quad \hbox{\rm and} \quad \ln \nu _1 \leq \alpha \ln \nu _2 ,
\]
from which comes
\[
\alpha =1 \quad \hbox{\rm and} \quad \ln \nu _2 = \ln \nu _1 \ ;
\]
thus the Lyapounoff exponents coincide on the periodic orbits.
\end{preuve}

The following result is a direct consequence of a property of transitive Anosov flows on~$3$-manifolds~\cite{L-M}.

\begin{cor}\label{cor:reguconj}
With the assumptions of the theorem~\ref{thm:reguconj}, the conjugacy~$h$  is of~${\cal C}^\infty $-class.
\end{cor}

The linearization allows to proof the following result.
\begin{thm}\label{thm:spm}
Let~$M$ be a closed (compact without boundary), connected, oriented surface.
Let~$g_1$ and~$g_2$ be two~${\cal C}^\infty $-riemannian metrics over~$M$
whose curvatures are negatively pinched~:~$-k_0^2\leq K_i\leq -k_1^2<0$ for~$i=1,2$,.
Let~$\kappa _1$ be a~${\cal C}^\infty $-magnetic field over~$M$.
The magnetic flow~$\psi ^1_t=\psi ^{g_1,\kappa _1}_t$ is supposed to have a negative Jacobi endomorphism 
(thus satisfying the pinching condition~\refeq{pince}).
If the magnetic flow~$\psi ^1_t$ and the geodesic flow~$\varphi ^2_t=\varphi ^{g_2}_t$ have the same marked length spectrum,
and if the surface~$M$ has the same total volume for the two metrics,
then the two metrics are isotopic, which means that one is the image of the other by a diffeomorphism~$f$ of~$M$,
homotopic to the identity, and the magnetic field~$\kappa _1$ is zero.
\end{thm}

\begin{preuve}
This results from the corollary~\ref{cor:reguconj}~: if the two volumes are equal, and the flows~$\psi ^1$ and~$\varphi ^2$ 
are conjugate by an absolutely continuous homeomorphism, the result is known~\cite{Gro1}.
\end{preuve}

To some extent, the linearisation determines the flow and the geometry.

\begin{thm}\label{thm:rigiR}
Let~$M$ be an oriented surface diffeomorphic to~$\R^2$, equipped with two~${\cal C}^\infty $-riemannian metrics~$g_1$ and~$g_2$
whose curvatures are negatively pinched~:~$-k_0^2\leq K_i\leq -k_1^2<0$ for~$i=1,2$. 
Let~$\kappa _1$,~$\kappa _2$ be two~${\cal C}^\infty $-magnetic fields over~$M$.
The two magnetic flows~$\psi ^1_t=\psi ^{g_1,\kappa _1}_t$ and~$\psi ^2_t=\psi ^{g_2,k_2}_t$ 
are supposed to have negative Jacobi endomorphisms satisfying the pinching condition~\refeq{pince}
(the~${\cal C}^1$-norms of~$\kappa _1$ and~$\kappa _2$ are thus bounded),
and the gradient of the centre-stable~$u_{-,i}$ spaces and the gradient of the centre-unstable~$u_{+,i}$ spaces for~$i=1,2$
are supposed to be uniformly bounded.
If there exist a diffeomorphism~$f:M\to M$ and a point~$p\in M$ satisfying
\[
\forall v \in T^1_{1,p}M \qquad
\exists v' \in T^1_{2,f(p)}M \qquad
E^2_{v'} \circ f = E^1_v ,
\]
then the two metrics are images one of each other by~$f$, and so are the two magnetic fields~:~$\kappa_2=\kappa_1 \circ f$.
\end{thm}
\begin{preuve}
Let~$w$ be in~$T^1_1M$,~$q=\pi w$ and~$v=v^1(p,w_{+\infty })$. Let~$v' $ be in~$T^1_{2,f(p)}M $ such that~$E^2_{v'} \circ f = E^1_v$.
We have necessarily~$d_wf(q)=v^2(f(q),v'_{+\infty })$~: thus the mapping~$f$ is an isometry.
Its differential conjugates the flows, therefore it preserves the geodesic curvature of the orbits,
which means the magnetic fields.
\end{preuve}

\begin{thm}\label{thm:rigicomp}
Let~$M$ be a closed (compact without boundary), connected, oriented surface.
Let~$g_1$ and~$g_2$ be two~${\cal C}^\infty $-riemannian metrics over~$M$
whose curvatures are negatively pinched~:~$-k_0^2\leq K_i\leq -k_1^2<0$ for~$i=1,2$.
Let~$\kappa _1$,~$\kappa _2$ be two~${\cal C}^\infty $-magnetic fields over~$M$.
The two magnetic flows~$\psi ^1_t=\psi ^{g_1,\kappa _1}_t$ and~$\psi ^2_t=\psi ^{g_2,\kappa _2}_t$ 
are supposed to have negative Jacobi endomorphisms (thus satisfying the pinching condition~\refeq{pince}).
If there exist two vectors~$v_1\in T^1_1\widetilde M$,~$v_2\in T^1_2\widetilde M$ 
and a~${\cal C}^1$-diffeomorphism~$f:M\to M$ homotopic to the identity,
of which a lift~$\widetilde f$ over~$\widetilde M$ satisfies~$E^2_{v_2}\circ \widetilde f = E^1_{v_1}$,
then the two metrics are isotopic, transported by~$f$,
and so are the two magnetic fields~:~$\kappa_2=\kappa_1 \circ f$.
\end{thm}
\begin{preuve}
Let~$w$ be in~$W^{CS}(v)$,~$q=\pi w$. We have necessarily~$d_w\widetilde f(q)=v^2(\widetilde f(q),v'_{+\infty })$. 
Since~$w$ is chosen arbitrarily in the centre-stable manifold of~$v$, 
we may chose a vector~$w$ whose projection on~$T^1_1M$ has a dense orbit under~$\psi ^1$ when the time tends to~$-\infty $~;
we may also replace~$w$ by every element of its orbit.
Writing down~$\Pi $ the covering of~$\widetilde M$ onto~$M$, we obtain~$d_{d\Pi (w)}f(\Pi q)=d\Pi ( v^2(\widetilde f(q),v'_{+\infty }) ) $.
The differential of~$f$ thus preserves the norm on a dense subset of the unitary tangent bundle~;
because it is assumed continuous, the mapping~$f$ is an isometry.
The differential of~$\widetilde f$ preserves the geodesic curvature on the whole orbit of~$w$,
thus by projecting and by a density argument, we deduce that~$\kappa_2$ is the composed of ~$\kappa_1$ by~$f$.
\end{preuve}

%%%%%%%%%%%%%%%%%%%%%%%%%%%%%%%%%%%%%%%%%%%%%%%%%%%

\end{document}